\documentclass{article}[12pt]
\usepackage{
latexsym,
amsthm,amsfonts,
epsfig,psfrag
}
\begin{document}


\newtheorem{lemma}{Lemma}
\newtheorem{theorem}{Theorem}
\renewcommand{\proofname}{Proof}

\def \H{{I\!\!H}}
\def \R{{I\!\!R}}
\newcommand{\conv}{\mathop{\rm conv}\nolimits}
\newcommand{\inter}{\mathop{\rm int}\nolimits}
\newcommand{\rel}{\mathop{\rm rel}\nolimits}
\newcommand{\ad}{\mathop{\rm ad}\nolimits}


\begin{center}
{\Large  \bf
Hyperbolic Coxeter $n$-polytopes with $n+3$ facets 
}
\vspace{17pt}

{\large P.~Tumarkin}

\vspace{19pt}

\parbox{10cm}{
{\bf Abstract.} A polytope is called a {\it Coxeter polytope} if its dihedral angles are
 integer parts of $\pi$. In this paper we prove that if a non-compact
 Coxeter polytope of finite volume in $\H^n$ has exactly $n+3$ facets
 then $n\le 16$. We also find an example in $\H^{16}$ and show that it
 is unique.  
}

\end{center}

\vspace{10pt}

{\bf 1.}
Consider a convex polytope $P$ in $n$-dimensional hyperbolic space
 $\H^n$. 

A polytope is called a {\it Coxeter polytope} if its dihedral angles are
 integer parts of $\pi$. Any Coxeter polytope $P$ is a fundamental domain
of the discrete group generated by the reflections with respect to
facets of $P$.

Of special interest are hyperbolic Coxeter polytopes of finite
volume. In contrast to spherical and Euclidean cases there is no
complete classification of such polytopes.
It is known that dimension of compact Coxeter polytope does not exceed $29$ 
(see \cite{V2}),  and dimension of non-compact polytope of finite
volume does not exceed $ 995$ (see \cite{995}).
Coxeter polytopes in
$\H^3$ are completely characterized by Andreev
\cite{A1},~\cite{A2}. There exists a complete  classification of
hyperbolic simplices~\cite{L},~\cite{V3} and hyperbolic $n$-polytopes
with $n+2$ facets~\cite{Kap} (see also~\cite{V1}),~\cite{Ess2},~\cite{arxiv}).

In~\cite{Ess2} Esselmann proved that  $n$-polytopes with $n+3$
facets do not exist in $\H^{n}$, $n>8$, and the example found by
Bugaenko~\cite{Bu2} in  $\H^{8}$ is unique. 
There is an example of finite volume non-compact polytope in
$\H^{15}$ with  $18$ facets (see~\cite{V1}).
The main result of this note is the following theorem:

\begin{theorem}
 Dimension of finite volume non-compact hyperbolic Coxeter
$n$-polytope with $n+3$ facets does not exceed $16$. In $\H^{16}$ there
exists a unique polytope with $19$ facets; its Coxeter diagram is
presented below.

\begin{figure}[htb!]
\begin{center}
\epsfig{file=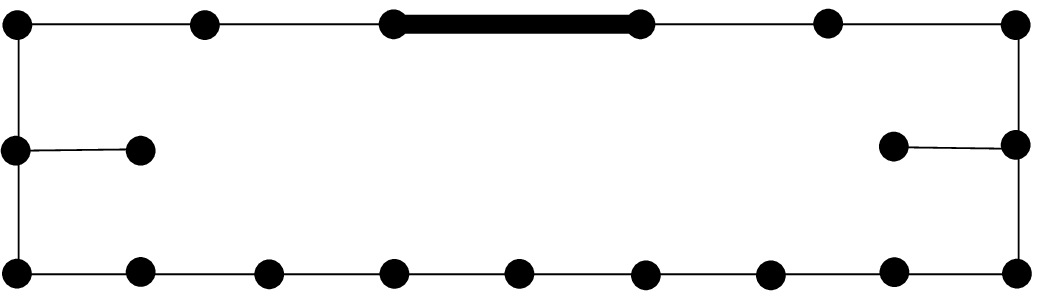,width=0.4\linewidth}
\end{center}
\label{16}
\end{figure}

\end{theorem}

\vspace{-15pt} 

{\bf 2.} 
To represent Coxeter polytopes one can use {\it Coxeter diagrams}.  
Nodes of Coxeter diagram correspond to facets of polytope. Two nodes
are joined by a $(m-2)$-fold edge or a $m$-labeled edge if the 
corresponding dihedral angle equals $\frac{\pi}{m}$. If the
corresponding facets are parallel the nodes are joined by a bold edge,
and if they diverge then the nodes are joined by a dotted edge
labeled by $\cosh(\rho)$, where $\rho$ is the distance between the facets.

Every combinatorial type of $n$-polytope with $n+3$ facets can be
represented by a standard two-dimensional {\it Gale diagram} (see, for
example,~\cite{gale}).  
This consists of vertices of regular $2k$-gon in $\R^{2}$ centered at
the origin and (possibly) the 
origin which are labeled according to the following rules:

1) Each label is a non-negative  integer, the sum of labels equals $n+3$.

2) Labels of neighboring vertices can not be equal to zero
   simultaneously.

3) Labels of opposite vertices can not be equal to zero
   simultaneously.

4) The points that lie in any open half-space bounded by a hyperplane
 through the origin have labels whose sum is at least two.

The combinatorial type of a convex polytope can be read off from 
the Gale diagram in the following way.
Each vertex $a_i$, $i=1,\dots,2k$, with label $\mu_i$ corresponds to
$\mu_i$ facets  
$f_{i,1},\dots,f_{i,\mu_i}$ of $P$. For any subset $I$ of the set of
facets of $P$ the intersection of facets 
$\{f_{j,\gamma} | (j,\gamma)\in I \}$ 
is a face of $P$ if and only if the origin is contained in the set 
$\conv \{a_{j} | (j,\gamma)\notin I\}$. 

By {\it pyramid} in $\H^n$ we mean a convex hull of a point (apex)
and an $(n-1)$-dimensional polytope that is not a simplex. It is easy
to see that a polytope $P$ is a pyramid if and only if the origin has
non-zero label in the  Gale diagram of $P$. In Section~3 we suppose
that $P$ is not a pyramid.   

{\bf 3.}
A connected Coxeter diagram $S$ is called a {\it Lann\'er (quasi-Lann\'er)
diagram} if  any subdiagram of $S$ is elliptic (elliptic or
parabolic), and  the diagram
$S$ is neither elliptic nor parabolic.
All Lann\'er diagrams are classified in~\cite{L}.

Let $G$ be the standard Gale diagram of polytope $P$. 
Denote by $S_{m,l}$ the following subdiagram of Coxeter diagram
$S(P)$:  $S_{m,l}$ corresponds to $l-m+1\ {\mbox {(mod
2k)}}$ consecutive vertices  $a_{m},\dots, a_{l}$ of $G$.

The following two lemmas can be easily derived from the definition of
Gale diagram and Theorems 3.1 and 3.2 of the paper~\cite{V1}.

\begin{lemma}\label{kv}
Suppose that the labels of vertices $a_i, a_{k+i}$ are not equal to
zero. Then\\
1) the labels of vertices $a_i$ and $a_{k+i}$ equal 1, and Coxeter
diagrams
$S_{i+1,k+i-1}$ and $S_{k+i+1,i-1}$ are connected and parabolic;\\
2) if the labels of vertices $a_{i+1}$ and $a_{k+i+1}$ are not equal
to zero then the Coxeter diagram $S_{i+1,k+i}$ is quasi-Lann\'er diagram;\\
3) if the label of vertex $a_{i+1}$ equals zero then the Coxeter diagram
$S_{i+2,k+i}$ is quasi-Lann\'er diagram. 

\end{lemma}

\begin{lemma}\label{l}
Suppose that labels of vertices $a_i, a_{k+i-1}$ are  equal to
zero. Then the Coxeter diagram $S_{i+1,k+i-2}$ is Lann\'er diagram.

\end{lemma}

Note that the number of vertices of any Lann\'er diagram does not
exceed $5$ (see \cite{L}), and the number of vertices of any quasi-Lann\'er diagram does not
exceed $10$ (see \cite{V3}). Using Lemma~\ref{l}  and statements 2) and
3) of Lemma~\ref{kv} we derive the following

\begin{lemma}\label{010}
Suppose that the labels of vertices $a_i$ and  $a_{i+2}$ are equal to
zero. Then the sum of labels of all vertices does not
exceed $20$.

\end{lemma}

An examination of the rest cases shows that dimension of polytope does
not exceed $17$. Thus, the sum of all labels of Gale diagram does
not exceed $20$, so the number of vertices with non-zero labels does
not exceed $20$, either.  From the other hand, if $k$ and the number of
vertices with  zero labels are big enough then the corresponding
Coxeter diagram contains at least two Lann\'er diagrams. In the latter
case $n\le 15$. This implies the restriction on $k$: if $n\ge 16$ 
then  $k\le 13$. Hence, we rest with a finite number of cases. Examining them
case by case we obtain the following result: there is no polytope for
$n=17$, and there is a unique one in $\H^{16}$.

{\bf 4.} Now we only need to check the case when  $P$ is a pyramid.

\begin{lemma}
Suppose that $P$ is a hyperbolic finite volume Coxeter $n$-polytope
with $n+3$ facets and $P$ is a pyramid. Then the combinatorial type of
$P$ is a pyramid over a product of three simplices.
 
\end{lemma}
\noindent
All the polytopes of this combinatorial type can be easily classified.
Their dimension does not exceed $11$.

\vspace{25pt}

\noindent
Independent Univ. of Moscow\\
e-mail:\qquad pasha@mccme.ru

\end{document}